\documentclass[11pt]{article}
\usepackage[dvips]{graphicx}
\usepackage{color,amsfonts,amssymb}
\usepackage{amsfonts,epsf,amsmath,tikz}
\usepackage{latexsym}
\usepackage{color}

\newtheorem{theorem}{Theorem}[section]
\newtheorem{lemma}[theorem]{Lemma}

\newcommand{\proof}{\noindent{\bf Proof.\ }}
\newcommand{\qed}{\hfill $\square$ \bigskip}

\newcommand{\cp}{\,\square\,}

\newcommand{\gp}{{\rm gp}}
\newcommand{\gpe}{{\rm gp_{e}}}

\newcommand{\diam}{{\rm diam}}
\newcommand{\ipe}{{\rm ip_e}}

\begin{document}
\title{Edge General Position Problem}

\author{
	Paul Manuel$^{a}$
	\and
	R.~Prabha$^{b}$
	\and
	Sandi Klav\v zar$^{c,d,e}$
}

\date{}

\maketitle
\vspace{-0.8 cm}
\begin{center}
	$^a$ Department of Information Science, College of Computing Science and Engineering, Kuwait University, Kuwait \\
	{\tt pauldmanuel@gmail.com}\\
	\medskip

	$^b$ Department of Mathematics, Ethiraj College for Women, Chennai, Tamilnadu, India \\
	{\tt prabha75@gmail.com}\\
	\medskip

	$^c$ Faculty of Mathematics and Physics, University of Ljubljana, Slovenia\\
	{\tt sandi.klavzar@fmf.uni-lj.si}\\
	\medskip

	$^d$ Faculty of Natural Sciences and Mathematics, University of Maribor, Slovenia\\
	\medskip

	$^e$ Institute of Mathematics, Physics and Mechanics, Ljubljana, Slovenia\\
	\medskip

\end{center}

\begin{abstract}
Given a graph $G$, the general position problem is to find a largest set $S$ of vertices of $G$ such that no three vertices of $S$ lie on a common geodesic. Such a set is called a ${\rm gp}$-$set$ of $G$ and its cardinality is the ${\rm gp}$-$number$, ${\rm gp}(G)$, of $G$. In this paper, the edge general position problem is introduced as the edge analogue of the general position problem.  The edge general position number, ${\rm gp_{e}}(G)$, is the size of a largest edge general position set of $G$. It is proved that ${\rm gp_{e}}(Q_r) = 2^r$ and that if $T$ is a tree, then ${\rm gp_{e}}(T)$ is the number of its leaves. The value of ${\rm gp_{e}}(P_r\cp P_s)$ is determined for every $r,s\ge 2$. To derive these results, the theory of partial cubes is used. Mulder's meta-conjecture on median graphs is also discussed along the way.  
\end{abstract}


\medskip
\noindent{\bf AMS Subj.\ Class.~(2020)}: 05C12, 05C76

%%%%%%%%%%%%%%%%%%%%%%%%%%%%%%%
\section{Introduction}
%%%%%%%%%%%%%%%%%%%%%%%%%%%%%%%

The geometric concept of points in general position, the still open Dudeney's no-three-in-line problem~\cite{Du17} from 1917 (see also~\cite{KuWo18, MiStSzSzZw16}), and the general position subset selection problem~\cite{FrKaNiNi17, PaWo13} from discrete geometry, all motivated the introduction of a related concept in graph theory as follows~\cite{MaKl18}. Let $G = (V (G), E(G))$ be a graph. Then the objective of the {\em general position problem} is to find a largest set of vertices $S \subseteq V(G)$, called a {\em $\gp$-set} of $G$, such that no vertex of $S$ lies on a geodesic between two other vertices of $S$. The {\em general position number} ($\gp$-{\em number} for short), $\gp(G)$, of $G$ is the cardinality of a gp-set of $G$. We point out that a couple of years earlier and in different terminology, the graph theory general position problem was considered in~\cite{ullas-2016} and that in the special case of hypercubes the problem has been much earlier studied by K\"orner~\cite{korner-1995}. 

Following the above listed seminal papers, the general position problem has been extensively investigated~\cite{anand-2019, ghorbani-2021, klavzar-2021+, klavzar-2021, klavzar-2019, MaKl18b, neethu-2021, patkos-2020, thomas-2020, tian-2021, tian-2020+}.  Let us emphasize some selected results. From~\cite{MaKl18} we recall that the general position problem is NP-complete and that in a block graph, the set of all simplicial vertices forms a gp-set. In~\cite{anand-2019} it is proved that $S \subseteq V(G)$ is a general position set if and only if the components of the subgraph induced by $S$ are complete subgraphs, the vertices of which form an in-transitive, distance-constant partition of $S$. In the same paper a formula for the gp-number of the complement of a bipartite graph is deduced and simplified for the complements of trees, of grids, and of hypercubes. In~\cite{klavzar-2019} the general position problem has been studied on different product graphs and connected with strong resolving graphs. It should be added that the concept of the general position set has been recently~\cite{klraya-2021+} extended to $d$-general position sets. 

The main topic of interest related to the general position problem thought was the Cartesian product operation. Let us denote by $X^n$ the  Cartesian product of $n$ factors each isomorphic to $X$, and let $P_\infty$ be the two-way infinite path. One of the main results from~\cite{MaKl18} asserts that $\gp(P_\infty^2) = 4$. In the same paper it was also proved that $10\le \gp(P_\infty^3) \le 16$. The lower bound $10$ was improved to $14$ in~\cite{klavzar-2021+}. These efforts culminated in~\cite{klavzar-2021} where it is proved that if $n\in {\mathbb N}$, then $\gp(P_\infty^n) = 2^{2^{n-1}}$. The general position problem in Cartesian products has been further investigated in~\cite{klavzar-2021+, tian-2021, tian-2020+}. In particular, it was proved in~\cite{tian-2021} that $\gp(G\cp H) \le n(G ) + n(H) - 2$, where the equality holds if and only if $G$ and $H$ are both generalized complete graphs. Moreover, the main result of~\cite{tian-2020+} asserts that the general position number is additive on Cartesian products of trees.

In this paper, the edge version of the graph theory general position problem is introduced. A set $S$ of edges of graph $G$ is said to be an {\em edge general position set} if no geodesic of $G$ contains three edges of $S$.  An edge general position set of maximum cardinality is called a $\gpe$-$set$ of the graph. An edge general position problem is to find a $\gpe$-set. The cardinality of a maximum edge general position set is called the \textit{edge general position number} (in short $\gpe$-number) of $G$, to be denoted by $\gpe(G)$. Our results indicate that the edge general position problem is a concept that deserves to be investigated, in particular, it is intrinsically different from the general position problem. 

We proceed as follows. In the next section we give further definitions needed and prove some preliminary results. In Section~\ref{sec:hyper} we first prove that $\gpe(Q_r) = 2^r$. This is in contrast with the fact that determining $\gp(Q_r)$ appears to be   very difficult~\cite{korner-1995}. We also prove that the leaves of a tree form an $\gpe$-set of it. These two results are then used to discuss Mulder's meta-conjecture on median graphs. In Section~\ref{sec:grids} we determine $\gpe(P_r\cp P_s)$ for all $r,s\ge 2$. Here the distinct difference between the vertex and the edge version of the general position problem is that the edge general position number of the $n \times n$ grid is proportional to $n$ whereas the general position number of the $n \times n$ grid remains constant even when $n$ tends to infinity. Moreover, we prove that the $\gpe$-set of $P_r\cp P_s$ is unique as soon as $r,s\ge 5$, another striking difference with the vertex version. 

%%%%%%%%%%%%%%%%%%%%%%%%%%%%%%%
\section{Preliminaries}
%%%%%%%%%%%%%%%%%%%%%%%%%%%%%%%

In this section we first state further concepts and the notation needed. We continue with some preliminary and auxiliary results. 

Unless stated otherwise, graphs considered in this paper are connected.  Let $G = (V(G), E(G))$ be a graph with vertex set $V(G)$ and edge set $E(G)$. Its order and size will be respectively denoted by $n(G)$ and $m(G)$. Let $P_n$ denote the path on $n$ vertices and $C_n$ the cycle on $n$ vertices. The {\em distance} $d_G(u, v)$ between vertices $u$ and $v$ of $G$ is the number of edges on a shortest $u, v$-path. Shortest paths are also known as {\em isometric paths} or {\em geodesics}. The {\em diameter} $\diam (G)$ of $G$ is the maximum distance between $u$ and $v$ of $G$. A subgraph $H$ of a graph $G$ is {\em isometric} if $d_H(x, y) = d_G(x, y)$ holds for every pair of vertices $x, y$ of $H$. A {\em pendant vertex} of a graph is a vertex of degree one, the edge incident to it is a {\em pendant edge}. 

The {\em Cartesian product} $G\cp H$ of graphs $G$ and $H$ is defined on the vertex set $V(G\cp H) = V(G)\times V(H)$, vertices $(g,h)$ and $(g',h')$ are adjacent if either $gg'\in E(G)$ and $h=h'$, or $g=g'$ and $hh'\in E(H)$. If $h\in V(H)$, then the subgraph of $G\cp H$ induced by the vertices $(g,h)$, $g\in V(G)$, is a {\em $G$-layer} and is denoted by $G^h$. Analogously, if $g\in V(G)$, then the {\em $H$-layer}\ $^g\!H$ is the subgraph of $G\cp H$ induced by the vertices $(g,h)$, $h\in V(H)$. $G$-layers and $H$-layers are isomorphic to $G$ and to $H$, respectively.

The {\em $r$-dimensional hypercube} $Q_r$, $r\ge 1$, is a graph with $V(Q_r) = \{0, 1\}^r$, and there is an edge between two vertices if and only if they differ  in exactly one coordinate. That is, if $x = (x_1, \ldots, x_r)$ and $y = (y_1, \ldots, y_r)$ are vertices of $Q_r$, then $xy\in E(Q_r)$ if and only if there exists  $j\in [r]$ such that $x_j \neq y_j$ and $x_i = y_i$ for every $i\ne j$. Note that $n(Q_r) = 2^r$ and $m(Q_r) = r2^{r - 1}$. Note also that $Q_r = Q_{r-1}\cp K_2$ holds for $r\ge 2$. 

If $\diam(G) = 2$, then a geodesic of $G$ contains at most two edges. Hence we have the following observation. 

\begin{lemma}
\label{lem:diam-2}
If $\diam(G) = 2$, then $\gpe(G) = m(G)$. 
\end{lemma}

Lemma~\ref{lem:diam-2} in particular implies that $\gpe(K_n) = \binom{n}{2}$ and that $\gpe(C_n) = n$ for $3\le n\le 5$. If $n \geq 6$, then it is easy to observe that $\gpe(C_n) = 4$.

An {\em isometric path edge cover} of a graph $G$ is a collection ${\cal P}$ of isometric paths of $G$ such that each edge of $G$ lies on at least one of the paths  from ${\cal P}$. The cardinality of a smallest isometric path edge cover is the {\em isometric path edge number} of $G$ and denoted by $\ipe(G)$. The following observation will turn out to be very useful, hence we state it as a lemma for further use. 

\begin{lemma}
\label{lem:ipe}
If $G$ is a connected graph, then $\gpe(G) \le 2\cdot\ipe(G)$.
\end{lemma}

\proof
Let ${\cal P}$ be an isometric path edge cover of $G$, where $|{\cal P}| = \ipe(G)$. Since paths from ${\cal P}$ are isometric, each of them contains at most two edges from an arbitrary edge general position set. Hence the conclusion. 
\qed

%%%%%%%%%%%%%%%%%%%%%%%%%%%%%%%
\section{Hypercubes, trees, and Mulder's meta-conjecture}
\label{sec:hyper}
%%%%%%%%%%%%%%%%%%%%%%%%%%%%%%%

A graph $G$ is a {\em partial cube} if $G$ is an isometric subgraph of some hypercube. Partial cubes have application in many different areas ranging from  interconnection networks~\cite{graham-1971}, media theory~\cite{eppstein-2008}, till mathematical chemistry, the papers~\cite{arockiaraj-2019, crepnjak-2017} are just a selection of many papers on the latter applications. For recent developments on the theory of partial cubes we refer to~\cite{chepoi-2020, marc-2017, mofidi-2020, polat-2020} and references therein. 

A key tool in the investigation of partial cubes is the Djokovi\'c-Winkler relation $\Theta$~\cite{djok-73,wink-84} defined as follows. Edges $xy$ and $uv$ of a graph $G$ are in relation $\Theta$ if $d_G(x,u) + d_G(y,v) \not= d_G(x,v) + d_G(y,u)$. Winkler~\cite{wink-84} proved that a connected graph $G$ is a partial cube if and only if $G$ is bipartite and $\Theta$ is transitive. As $\Theta$ is reflexive and symmetric and the edge set of an arbitrary graph, it partitions the edge set of a partial cube into $\Theta$-{\em classes}. 

\begin{lemma}
\label{lem:two-classes}
Let $G$ be a partial cube and let $F_1$ and $F_2$ be $\Theta$-classes of $G$. Then $F_1 \cup F_2$ is an edge general position set of $G$. 
\end{lemma}

\proof
It is well-known that if $P$ is a shortest path in a graph $G$, then no two edges of $P$ are in relation $\Theta$, cf.~\cite[Lemma 11.1]{HIK-2011}. Let $e, f, g\in F_1 \cup F_2$. Suppose first that $e, f, g\in F_1$. Then no two of these edges lie on a common geodesic. The case $e, f, g\in F_2$ is analogous. Suppose second that, without loss of generality, $e, f\in F_1$ and $g \in F_2$. But then $e$ and $f$ are not on a common geodesic. In any case, $e$, $f$, and $g$ are not on a common geodesic.  
\qed

\begin{figure}[ht!]
\begin{center}
\begin{tikzpicture}[scale=1.0,style=thick]
\def\vr{3pt}
%% vertices defined %%
\path (0,0) coordinate (x1);
\path (0,1) coordinate (x2);
\path (0,2) coordinate (x3);
\path (0,3) coordinate (x4);
\path (1,0) coordinate (x5);
\path (1,1) coordinate (x6);
\path (1,2) coordinate (x7);
\path (1,3) coordinate (x8);
\path (2,0) coordinate (x9);
\path (2,1) coordinate (x10);
\path (3,0) coordinate (x11);
\path (3,1) coordinate (x12);
%% edges %%
\draw (x1) -- (x2) --(x3) -- (x4) -- (x8) -- (x7) -- (x6) -- (x10) -- (x12) -- (x11) --(x9) -- (x5) -- (x1);
\draw (x3) -- (x7);
\draw (x2) -- (x6);
\draw (x5) -- (x6);
\draw (x9) -- (x10);
\draw (0.5,-0.15) -- (0.5,0.15);
\draw (0.5,0.85) -- (0.5,1.15);
\draw (0.5,1.85) -- (0.5,2.15);
\draw (0.5,2.85) -- (0.5,3.15);
\draw (-0.15,0.55) -- (0.15,0.55); \draw (-0.15,0.45) -- (0.15,0.45); 
\draw (0.85,0.55) -- (1.15,0.55); \draw (0.85,0.45) -- (1.15,0.45); 
\draw (1.85,0.55) -- (2.15,0.55); \draw (1.85,0.45) -- (2.15,0.45); 
\draw (2.85,0.55) -- (3.15,0.55); \draw (2.85,0.45) -- (3.15,0.45); 
%% vertices %%
\draw (x1)  [fill=white] circle (\vr);
\draw (x2)  [fill=white] circle (\vr);
\draw (x3)  [fill=white] circle (\vr);
\draw (x4)  [fill=white] circle (\vr);
\draw (x5)  [fill=white] circle (\vr);
\draw (x6)  [fill=white] circle (\vr);
\draw (x7)  [fill=white] circle (\vr);
\draw (x8)  [fill=white] circle (\vr);
\draw (x9)  [fill=white] circle (\vr);
\draw (x10)  [fill=white] circle (\vr);
\draw (x11)  [fill=white] circle (\vr);
\draw (x12)  [fill=white] circle (\vr);
%% others %%
\end{tikzpicture}
\end{center}
\caption{A partial cube $G$.}
\label{fig:pc}
\end{figure}
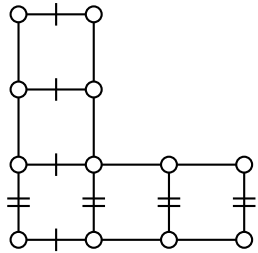

As a small example consider the partial cube $G$ from Fig.~\ref{fig:pc}. It has two $\Theta$-classes each containing four edges (marked on the figure), hence Lemma~\ref{lem:two-classes} implies $\gpe(G) \ge 8$. On the other hand, it is not difficult to find an isometric path edge cover of $G$ consisting of four geodesics, hence $\gpe(G) \le 8$ by Lemma~\ref{lem:ipe}. We conclude that $\gpe(G) = 8$.

\begin{theorem} 
\label{thm:hypecubes}
If $r\ge 2$, then $\gpe(Q_r) = 2^r$.
\end{theorem}

\proof 
It is well-known that hypercubes are partial cubes. Moreover, a $\Theta$-class of a hypercube is formed by the edges whose endpoints differ in the same, fixed coordinate. Hence, $Q_r$ has $r$ $\Theta$-classes, each containing $2^{r-1}$ edges. Lemma~\ref{lem:two-classes} thus implies that $\gpe(Q_r) \ge 2^r$. 

\begin{figure}[ht!]
\begin{center}
\begin{tikzpicture}[scale=1.0,style=thick]
\def\vr{3pt}
%% vertices defined %%
\path (0,0) coordinate (x1);
\path (0,4) coordinate (x2);
\path (1,1) coordinate (x3);
\path (1,3) coordinate (x4);
\path (3,1) coordinate (x5);
\path (3,3) coordinate (x6);
\path (4,0) coordinate (x7);
\path (4,4) coordinate (x8);
%% edges %%
\draw[thick,->] (0,4) -- (0,0.1);
\draw[thick,->] (0,0) -- (0.9,0.9);
\draw[thick,->] (1,1) -- (2.9,1);
\draw[thick,dotted, ->] (1,1) -- (1,2.85);
\draw[thick,dotted, ->] (1,3) -- (0.1,3.9);
\draw[thick,dotted, ->] (0,4) -- (3.9,4);
\draw[thick,densely dotted, ->] (3,3) -- (3,1.1);
\draw[thick,densely dotted, ->] (3,1) -- (3.9,0.1);
\draw[thick,densely dotted, ->] (4,0) -- (0.1,0);
\draw[thick,dashed, ->] (4,0) -- (4,3.9);
\draw[thick,dashed, ->] (4,4) -- (3.1,3.1);
\draw[thick,dashed, ->] (3,3) -- (1.1,3);
%% vertices %%
\draw (x1)  [fill=white] circle (\vr);
\draw (x2)  [fill=gray] circle (\vr);
\draw (x3)  [fill=gray] circle (\vr);
\draw (x4)  [fill=white] circle (\vr);
\draw (x5)  [fill=white] circle (\vr);
\draw (x6)  [fill=gray] circle (\vr);
\draw (x7)  [fill=gray] circle (\vr);
\draw (x8)  [fill=white] circle (\vr);
%% others %%
\draw (0.4,2) ellipse (1cm and 3cm);
\draw (3.6,2) ellipse (1cm and 3cm);

% \draw [left] (y1) node {$y_1$};
\end{tikzpicture}
\end{center}
\caption{Covering $Q_3$ with four shortest path. The paths are further oriented such that the endpoints of the paths form one of the bipartition sets.}
\label{fig:Q3}
\end{figure}
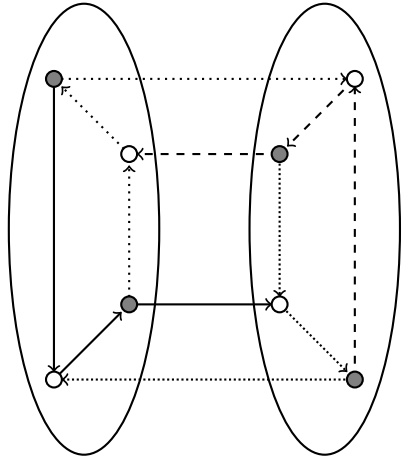

To prove the reverse inequality we claim the following: if $r\ge 2$, then there exists an isometric path edge cover ${\cal P}_r$ of $Q_r$, where the paths from ${\cal P}_r$ can be oriented such that the endpoints of the paths form one of the bipartition sets of $Q_r$. It is straightforward to check that the claim holds for $r=2$, see the left hand side of Fig.~\ref{fig:Q3} where the circled $Q_2$ is covered with two $2$-paths. These two paths end in the two gray vertices which form a bipartition set of this $Q_2$. 

Let now ${\cal P}'_r$ be an isometric path edge cover of $Q_r$, $r\ge 2$, together with the orientation of the paths such that the endpoints of the paths form one of the bipartition sets of $Q_r$. Let further ${\cal P}_r''$ be an isometric edge path cover of $Q_r$, where the paths from ${\cal P}_r''$ are orientated such that the endpoints of the paths form the other bipartition set of $Q_r$. The existence of ${\cal P}'_r$ will be guaranteed by the induction, while ${\cal P}''_r$ can be constructed from ${\cal P}'_r$ by applying the automorphism of $Q_r$ which assigns to each vertex $u$ the vertex $u + 00\ldots 01$, where $+$ stands for the component-vise summation modulo $2$. (See the right-hand side of Fig.~\ref{fig:Q3}, where ${\cal P}_2''$ is shown with the endpoints of its two shortest paths drawn gray again.) Consider now $Q_{r+1}$ as the Cartesian product $Q_r\cp K_2$ with layers $Q_r^1$ and $Q_r^2$, and respective isometric path edge covers ${\cal P}_r'$ and ${\cal P}_r''$. Then extend each path from ${\cal P}_r'$ by the edge from its endpoint in $Q_r^1$ to its (unique) neighbor in  $Q_r^2$ and orient the new edge in this direction. Note that the new path is a shortest path of $Q_{r+1}$. Similarly, extend each path from ${\cal P}_r''$ by the edge from its endpoint in $Q_r^2$ to its (unique) neighbor in  $Q_r^1$ and orient the new edge in this direction. Note that the ends of the paths obtained by extending paths from ${\cal P}_r'$ together with the ends of the paths obtained by extending paths from ${\cal P}_r''$ form a bipartition set of $Q_{r+1}$. (See Fig.~\ref{fig:Q3} again for these extensions and note that their ends form the bipartition set of $Q_3$ drawn white.) This proves the claim. 

It follows from the above proved claim that $Q_r$ admits an isometric path edge cover ${\cal P}_r$ together with an orientation which reveals that $|{\cal P}_r|$ is the size of a bipartition set of $Q_r$, that is, $|{\cal P}_r| = 2^{r-1}$. By Lemma~\ref{lem:ipe} we conclude that $\gpe(Q_r) \le 2^r$. 
\qed

Another fundamental class of partial cubes is the class of trees for which we have the following. 

\begin{theorem} 
\label{thm:trees}
If $L$ is the set of pendant edges of a tree $T$, then $\gpe(T) = |L|$.
\end{theorem}

\proof 
We first prove that $\gpe(T) \le |L|$. Let $S$ be an edge general position set of $T$ and suppose that a non-pendant edge, say $e$, belongs to $S$. Let $T_1$ and $T_2$ be the components of $T - e$. Then  $E(T_1) \cap S = \emptyset$ or $E(T_2) \cap S = \emptyset$. Indeed, if there would be edges $e'\in E(T_1) \cap S$ and $e''\in E(T_2) \cap S$, then the edges $e$, $e'$, and $e''$ would lie on a common geodesic. We may thus without loss of generality assume that $E(T_2 ) \cap S = \emptyset$. Let $f$ be a pendant edge of $T_2$ that is at largest possible distance from $e$. (Here an arbitrary pendant edge of $T_2$ would actually do the job, except maybe a pendant edge adjacent to $e$.) Then it is straightforward to see that $(S \setminus \{e\}) \cup \{f\}$ is an edge general position set of $T$. Inductively continuing this process we end up with an edge general position set of $T$ which contains only pendant edges and has the same cardinality as $S$. We conclude that $\gpe(T) \le |L|$.

To see that $\gpe(T) \ge  |L|$ holds, observe that the set of pendant edges of a tree form an edge general position set because a geodesic of $T$ can pass through at most two pendant edges. 
\qed
	
Trees and median graphs are fundamental building blocks of the class of median graphs, cf.~\cite{klavzar-1999}, which is in turn (probably) the most important subclass of partial cubes. In 1990, Mulder proposed the following meta-conjecture: Any (sensible) property that is shared by trees and hypercubes is shared by all median graphs, see~\cite{mulder-2016}. Theorems~\ref{thm:hypecubes} and~\ref{thm:trees} do not share some obvious common property, that is, for hypercubes $\gpe$-sets constructed are the union of two $\Theta$-classes, while for trees their $\gpe$-sets are the sets of their leaves. Hence it is yet to be seen whether Mulder's meta-conjecture applies to the edge general position number of median graphs. We will further comment on this at the end of the next section.    

%%%%%%%%%%%%%%%%%%%%%%%%%%%%%%%
\section{Grid networks}
\label{sec:grids}
%%%%%%%%%%%%%%%%%%%%%%%%%%%%%%%

In this section we determine the $\gpe$-number of Cartesian products of two paths, known also as grid networks. We will always assume that $V(P_r) = [r] = \{1,\ldots, r\}$. The $P_r$-layers of $P_r \cp P_s$ are thus denoted by $P_r^i$, $i\in [s]$, and the $P_s$-layers by $^j\!P_s$, $j\in [r]$.  An edge $e = uv$ of $P_r \cp P_s$ is called a {\em boundary edge} if $\{d(u), d(v)\} = \{2, 3\}$ or $\{d(u), d(v)\} = \{3\}$, it is a {\em semi-boundary edge} if $\{d(u), d(v)\} = \{3, 4\}$, otherwise $e$ is an {\em internal edge}; that is, $e$ is internal if $\{d(u), d(v)\} = \{4\}$.

The main result of this section reads as follows. By the commutativity of the Cartesian product operation, the result covers all non-trivial grid networks, that is, we may without loss of generality assume that $r\ge s$.

\begin{theorem}
\label{thm:grids}
If $r\ge s\ge 2$, then 
$$\gpe(P_r \cp P_s) = 
\begin{cases}
r + \normalfont{2}; & s = 2,\\
2r; & s = 3,\\
2r + 2s - 8; & s \geq 4\,.
\end{cases}
$$
\end{theorem}

\begin{figure}[ht!]
\begin{center}
\begin{tikzpicture}[scale=1.0,style=thick]
\def\vr{3pt}
%% vertices defined %%
\foreach \j in {0,1,...,6}
\foreach \i in {0,1,...,9}
{
  \path (\i,\j) coordinate (x\i\j);
}
% \path (0,0) coordinate (x1);
%% edges %%
\foreach \j in {0,1,...,6}
{
   \draw (x0\j) -- (x9\j);
}   
\foreach \i in {0,1,...,9}
{
   \draw (x\i0) -- (x\i6);
}   
%% vertices %%
\foreach \j in {0,1,...,6}
\foreach \i in {0,1,...,9}
{
   \draw (x\i\j)  [fill=white] circle (\vr);
}
%% paths of the cover
\draw[thick] (2.15,0.1) -- (2.15,6.15);
\draw[thick] (0.1,6.15) -- (2.15,6.15);
\draw[thick, dashed] (0.85,0.1) -- (0.85,6.25);
\draw[thick, dashed] (0.85,6.25) -- (9,6.25);
\draw[thick] (6.85,-0.15) -- (6.85,5.85);
\draw[thick] (6.85,-0.15) -- (8.85,-0.15);
\draw[thick, densely dotted] (8.15,5.85) -- (8.15,-0.25);
\draw[thick, densely dotted] (8.15,-0.25) -- (0,-0.25);
\end{tikzpicture}
\end{center}
\caption{Part of an isometric path edge cover of $P_{10}\cp P_7$.}
\label{fig:grid}
\end{figure}
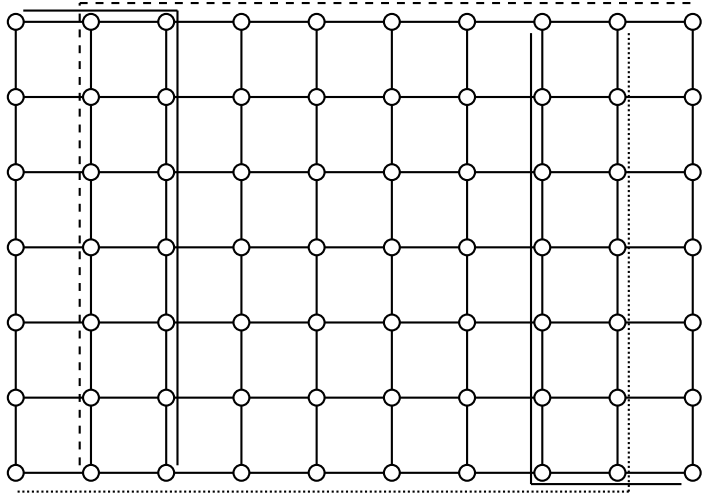

\proof
Suppose first that $s\ge 6$. Construct an isometric path edge cover of $P_r\cp P_s$ as follows. Select first four paths that cover all the edges from the $P_s$-layers $^2\!P_s$, $^3\!P_s$, $^{r-2}\!P_s$, $^{r-1}\!P_s$, as well as all the edges from the $P_r$-layers $P_r^1$ and $P_r^s$. See Fig.~\ref{fig:grid} where such four path are drawn for the case $r=10$ and $s=7$. 

\noindent
By symmetry, select additional four paths that cover all the edges from the $P_r$-layers $P_r^2$, $P_r^{3}$, $P_r^{s-2}$, $P_r^{s-1}$, as well as all the edges from the $P_s$-layers $^1P_s$ and $^rP_s$. Note that the eight paths selected so far cover all the edges from six $P_r$-layers and all the edges from six $P_s$-layers. Hence we may easily complete the path edge cover by adding $r-6$ paths that cover the edges not yet covered in the $P_s$-layers and $s-6$ paths that cover the edges not yet covered in the $P_r$-layers. The constructed path edge cover contains $8 + (r-6) + (s-6) = r+s-4$ paths. Lemma~\ref{lem:ipe} implies that $\gpe(P_r\cp P_s) \le 2r + 2s - 8$. On the other hand, the set of the semi-boundary edges of $P_r\cp P_s$ is an edge general position set of cardinality $2r+2s-8$. We conclude that it is a $\gpe$-set, that is, $\gpe(P_r\cp P_s) = 2r + 2s - 8$ when $r\ge s\ge 6$. 

Suppose next that $s=5$. Let $F$ be an arbitrary edge general position set of $P_r \cp P_5$ and distinguish two cases. 

\medskip\noindent
{\bf Case 1}. $|F\cap E(^iP_5)| \le 1$ for each $i\in [r]$.\\
Suppose first that $|F\cap E(\bigcup_i{^iP_5})| = r-x$, where $x\ge 3$. This means that there are exactly $x$ $P_s$-layers with no edges from $F$ and $r-x$ $P_s$-layers with exactly edges from $F$. Since every of the five $P_r$-layers contains at most two edges from $F$, it follows that $|F| \le (r-x) + 10$. We wish to show that $(r-x) + 10 \le 2r+2$ which is equivalent that $r+x\ge 8$ holds. Since $r\ge 5$ and $x\ge 3$, this is indeed the case. 

Suppose second that $|F\cap E(\bigcup_i{^iP_5})| = r-x$, where $x\in \{0,1,2\}$. If each $P_r$-layer contains at most one edge from $F$, then $|F| \le (r-x) + 5$. We wish to show that $r-x+5 \le 2r + 2$, which is equivalent to $r+x\ge 3$, the latter being clearly true. Hence suppose that at least one $P_r$-layer shares two edges with $F$. In this case we infer that there are at least two $P_5$-layers that have no edge in $F$. If follows that $|F| \le (r-2) + 10$ and we wish that $(r-2) + 10\le 2r + 2$, that is, $r\ge 6$. This settles Case 1 in all possibilities except when $r=5$. 

It remains to consider $P_5\cp P_5$. The theorem asserts that $\gpe(P_5\cp P_5) = 12$. The $12$ semi-boundary edges imply that $\gpe(P_5\cp P_5) \ge 12$. Here it remains to prove that $\gpe(P_5\cp P_5) \le 12$ under the assumption of Case 1. This is clearly the case if each $P_r = P_5$-layer contain at most one edge from $F$. Hence assume that some $P_r = P_5$-layer contains two edges from $F$.  Then we see similarly as above that the $P_s$-layers ($s=5$ of course) contain at most $3$ edges from $F$. Hence the only possibility to get more than $13$ edges in $F$ is that each $P_r = P_5$-layer contains exactly two edges from $F$. In this case we first infer that these $10$ edges must project to exactly two edges of the $P_r = P_5$ factor. But then it is clear that there are no more edges in $F$. 

\medskip\noindent
{\bf Case 2}. There exists an $i\in [r]$ such that $|F\cap E(^iP_5)| = 2$.\\
Suppose first that this $i$ is unique. Since this $P_5$-layer has two edges from $F$, we see that at least two of the $P_r$-layers have no edges in $F$. Since the other have at most two such edges, we conclude that $|F| \le 2 + (r-1) + 3\cdot 2$ holds in this subcase. We wish to show that $2 + (r-1) + 6 \le 2r + 2$, which is equivalent to $r\ge 5$. So we are done in this subcase. 

Suppose second that there are at least two $P_5$-layers with exactly two edges from $F$. Then we infer again that the projection of these edges on the factor $P_5$ contain exactly two edges. But then we easily obtain that in each case $|F|\le 2r + 2$. For instance, if each of the $P_5$ layers contain exactly two edges from $F$, then actually none of the edges from the $P_r$-layers lie in $F$. This settles Case 2 and hence theorem is proved for the case $s=5$. 

\medskip
Suppose next that $s=4$, that is, $P_s = P_4$. Let $F$ be an arbitrary edge general position set of $P_r \cp P_4$. 

\medskip\noindent
{\bf Case 1}. $|F \cap E(^iP_4)|\le 1$ for each $i \in [r]$. 

Suppose first that $|F \cap E(\bigcup_i {^iP_4}) | = r-x$ where $x \ge 4$. This means that there are exactly $x$ $P_s$-layers with no edges from $F$. Since each of the four $P_r$-layers contains at most two edges from $F$, it follows that $|F| \le (r-x)+8$. We wish to show that $(r-x)+8 \le 2r$ which is equivalent to $r+x\ge 8$. Since $r\ge 4$ and $x\ge 4$ this is indeed the case. 

Suppose second that $F\cap E(\bigcup_i {^iP_5}) = r-x$ where $x \in \{0, 1, 2, 3\}$. If each $P_r$-layer contains at most one edge from $F$, then $|F|\le (r-x)+4$. We wish to show that $r-x+4\le 2r$ which is equivalent to $r+x\ge 4$, the latter being clearly true. Hence suppose that at least one $P_r$-layer shares two edges with $F$. In this case we infer that there are at least two $P_4$-layers that have no edge in $F$. It follows that $|F| \le (r-2)+8$. Since we would like to see that $(r-2)+8\le 2r$,  that is $r\ge 6$, this settles Case 1 in all possibilities except when $r\in \{4,5\}$. 

Consider $P_4\cp P_4$. The theorem asserts that $\gpe(P_4\cp P_4)=8$ and the eight semi-boundary edges imply that $\gpe(P_4\cp P_4)\ge 8$. Hence it remains to prove that $\gpe(P_4\cp P_4)\le \ 8$ (under the assumption of Case 1). This is clearly the case if each $P_r = P_4$-layer contains at most one edge from $F$. Hence assume that some $P_r = P_4$-layer contains two edges from $F$. Then there exist at least two $P_s = P_4$-layers that have no edge in $F$. Hence, under the assumption of Case 1, the $P_s$-layers together contain at most two edges from $F$. Hence the only possibility to get more than eight edges in $F$ is that at least three among the $P_r = P_4$-layers must contain exactly two edges from $F$. But then the $P_s$-layers can contain at most one edge from $F$ and we easily conclude that $|F|\le 8$. 

Consider $P_5\cp P_4$. The theorem asserts that $\gpe(P_5\cp P_4) = 10$. The ten semi-boundary edges imply that $\gpe(P_5\cp P_4)\geq 10$. It remains to prove that $\gpe(P_5\cp P_4)\leq 10$ under the assumption of Case 1. If each $P_r = P_5$-layer contains at most one edge from $F$, then there is nothing to prove. Hence assume that some $P_r = P_5$-layer contains two edges from $F$. Then the $P_s = P_4$-layers together contain at most three edges from $F$. Hence the only possibility to get more than ten edges in $F$ is that at all four $P_5$-layers contain exactly two edges from $F$. But in that case, no $P_4$-layer can share an edge with $F$. 

\medskip\noindent
{\bf Case 2}. There exists an $i \in [r]$ such that $|F \cap E(^iP_4)| =2$.

Suppose first that this $i$ is unique. Since this $P_4$-layer has two edges from $F$, we see that at least two of the $P_r$-layers have no edges in $F$. Since the other $P_r$-layers have at most two such edges, we conclude that $|F| \le 2+(r-1) +2\cdot 2$ holds in the subcase. We wish to show that $2 + (r-1) +4 \le 2r$ which is equivalent to $r\ge 5$. This settles Case 2 in all possibilities except when $r= 4$. For $P_4\cp P_4$ we can use a parallel argument as we gave in Case 1. 

Suppose second that there are at least two $P_4$-layers with exactly two edges from $F$. Then we infer that the projection of these edges on the factor $P_4$ contains exactly two edges. But then we can verify easily that in each case, $|F|\le 2r$. For instance, if each of the $P_s = P_4$-layers contain exactly two edges from $F$, then actually none of the edges from the $P_r$-layers lie in $F$. This settles Case 2 and hence the theorem for $s=4$.

\medskip
Next, let $s=3$. Note that the theorem asserts that $\gpe(P_r\cp P_4) = \gpe(P_r\cp P_3)  = 2r$. As we have already proved that $\gpe(P_r\cp P_4) = 2r$ and $P_r\cp P_3$ is an isometric subgraph of $P_r\cp P_4$, it follows that $\gpe(P_r\cp P_3) \le \gpe(P_r\cp P_4)$. On the other hand, a set of $2r$ edges that is a $\gpe$-set of $P_r\cp P_4$ can also be used as an edge general position set of $P_r\cp P_3$, that is, a $\gpe$-set of $P_r\cp P_3$.

\medskip
Finally, let $s=2$. If each of the two $P_r$-layers intersects $F$ in at most one edge, then clearly $|F| \le r + 2$. On the other hand, if one of the two $P_r$-layers contains two edges from $F$, then at least two edges of the $P_2$-layers are not in $F$. Hence again $|F| \le 2 + 2 + (r-2) = r + 2$. On the other hand, the union of two $\Theta$-classes, one with $r$ edges, and the other with $2$ edges, is an edge general position set by Lemma~\ref{lem:two-classes}. 
\qed

We next supplement Theorem~\ref{thm:grids} by the following information. 

\begin{theorem} 
	\label{thm2 grid}
	If $r,s\ge 5$, then the $\gpe$-set of $P_r \cp P_s$ is unique.
\end{theorem}

\proof 
By Theorem~\ref{thm:grids} we know that the set $S$ of semi-boundary edges of $P_r \cp P_s$ form a $\gpe$-set of $P_r \cp P_s$. We need to prove that there is no other $\gpe$-set. For this sake assume that $T$ is an arbitrary edge general position set different from $S$. Our goal is to show that $|T| < 2r + 2s - 8$. Since $T\ne S$ we see that $T$ contains a boundary or an internal edge and distinguish our considerations accordingly.  
		
Suppose first that $T$ contains a boundary edge $e$. By the commutativity of $P_r\cp P_s$ and by the symmetry between the layers  $^1\!P_s$ and $^r\!P_s$ we may without loss of generality assume that $e \in E(^1\!P_s)$. 

\medskip\noindent
{\bf Case 1}. $T\cap E(P_r^1) \neq \emptyset$ or $T\cap E(P_r^s) \neq \emptyset$. \\
In this case we infer that $|T\cap E(P_r^i)| \leq 1$ for $i \in [s]$ and that $|T\cap E(^j\!P_s)| \leq 1$ for $j\in [r]$. Then clearly, $|T| \leq r + s$ and we wish to show that $r + s < 2r + 2s - 8$. We are done because $r + s > 8$ holds. 

In the rest of the argument we may hence assume that $T\cap E(P_r^1) = \emptyset$ and $T\cap E(P_r^s) = \emptyset$. 

\medskip\noindent
{\bf Case 2}. $|T \cap E(^1\!P_s)| = 1$ and $|T\cap E(^r\!P_s)| \in \{0,1\}$. 
In this case $|T\cap E(P_r^i)| \leq 1$ holds for every $2\leq i \leq s-1$. This in turn implies that $|T| \leq 2(r - 2) + 2 + s - 2 = 2r + s - 4$. We wish to show that $2r + s - 4 < 2r +2s - 8$, and this indeed holds since $s > 4$. 

\medskip\noindent
{\bf Case 3}. $|T \cap E(^1\!P_s)| = 2$ and $T\cap E(^r\!P_s) = \emptyset$. \\
Then clearly  $T\cap E(P_r^1) = T\cap E(P_r^s) = \emptyset$ which implies $|T| \leq 2(r - 1) + s - 2$. We wish that $2(r - 1) + s - 2 < 2r + 2s - 8$ which implies $s > 4$ which is indeed the case.

\medskip\noindent
{\bf Case 4}. $|T \cap E(^1\!P_s)| = 2$ and $|T\cap E(^r\!P_s)| \in \{1,2\}$. \\
In this case there exists at least three $P_r$-layers that contribute no edges to $T$. Hence $|T| \leq 2(r - 1) + 1 + s - 3 = 2r +s - 4$. Now $2r + s - 4 < 2r +2s - 8$ implies $s > 4$ which is indeed true.

\medskip
Suppose second that $T$ contains an internal edge. We may further assume that $T$ contains no boundary edge as we have already dealt with this situation. Without loss of generality suppose there exists a $P_s$-layer, say $^k\!P_s$, such that $T \cap E(^k\!P_s) = \{e_1, e_2\}$, where at least one of $e_1$ and $e_2$ is an internal edge. Now observe that at least three $P_r$-layers (including two layers $P_r^1$ and $P_r^s$) contribute no edges to $T$. Hence $|T| \leq 2(r - 2) + 2(s - 3) < |S|$.
\qed

Note that Theorem~\ref{thm2 grid} does not hold for the grids $P_4 \cp P_4$ and $P_5 \cp P_4$. Indeed, the sets of semi-boundary edges of then contain $8$ and $10$ edges, respectively. On the other hand, each $\Theta$-class of $P_4 \cp P_4$ has four edges, so the union of arbitrary two $\Theta$-classes is also a $\gpe$-set of $P_4 \cp P_4$.  In  $P_5 \cp P_4$, there are three $\Theta$-classes with five edges each, hence also the union of any two of them is a $\gpe$-set of $P_5 \cp P_4$. 

Let us return to the Mulder's meta-conjecture. In the previous section we found out that unions of two $\Theta$-classes are $\gpe$-sets of hypercubes and that the unique $\gpe$-set of a tree is the set of its pendant edges. Since each edge of a tree forms a $\Theta$-class, one could say that a (weak) common property of hypercubes and trees is that $\gpe$-sets are unions of $\Theta$-classes. But this does not extend to all median graphs. Since the Cartesian product of two median graphs is again median, we see that $P_r\cp P_s$ is a median graph. However, we have just seen that in $P_r\cp P_s$, where $r\ge s\ge 5$, the set of its semi-boundary edges forms a unique $\gpe$-set. Clearly, this set is not a union of $\Theta$-classes of $P_r\cp P_s$.

%%%%%%%%%%%%%%%%%%%%%%%%%%%%%%%
\section{Conclusions}
%%%%%%%%%%%%%%%%%%%%%%%%%%%%%%%

In this paper, the edge general position problem is completely solved for hypercubes and two-dimensional grids. A notable contribution of this paper is the discussion of Mulder's meta-conjecture on median graphs from the perspectives of the edge general position problem. This problem may be studied for other classes of graphs such as Cayley graphs, perfect graphs, bipartite graphs, etc.

It would also be pertinent to view this problem in generalized perspectives. For instance, for a given integer $k\ge 3$, one may call a set $S$ of edges of a graph $G$ an {\em edge $k$-general position set} if no $k$ edges of $S$ lie on a common geodesic. An edge $k$-general position set $S$ of maximum cardinality is an  $k$-$\gpe$-{\em set} of $G$ and its cardinality is the {\em edge $k$-general position number} (in short, $k$-$\gpe$-{\em number}) of $G$ and is denoted by $k$-$\gpe(G)$. When $k = 3$, this problem becomes an edge general position problem. The complexity status of edge $k$-general position problems is not known. Enthusiastic graph theory students will find the ``edge $k$-general position problem" an interesting topic for further research. 

\section*{Acknowledgements}

The research is supported by Kuwait University, Kuwait.

\end{document}